\documentclass[a4paper,10pt]{article}
\usepackage[top=1in, bottom=1in, left=1in, right=1in]{geometry}
\usepackage{amsfonts}
\usepackage{amssymb}
\usepackage{amsthm}
\usepackage{amsmath,amssymb,latexsym} 
\usepackage{fancyhdr}
\usepackage{graphicx, graphics}

\newtheorem{thm}{Theorem}
\newtheorem{deff}{Definition}
\newtheorem{prop}[thm]{Proposition}

\newtheorem{cor}[thm]{Corollary}

\begin{document}

\title{Notes on algebraic structure of truth tables of bracketed formulae connected by implications }
\author{\textit{Volkan Yildiz} }
\date{  }
\maketitle

\abstract
In this paper we investigate the algebraic structure of the truth tables
of all bracketed formulae with $n$ distinct variables connected by the binary connective
of implication.
\\\\\\

Keywords: Truth Tables, Implication, Catalan numbers, Generating Functions,  Monoids\\
AMS classification: 05A15, 05A16, 03B05, 11B75, 20M10

\pagebreak

\begin{flushright}{ \it{for Ares}}\end{flushright}

\section{Introduction}
We explore the algebraic structure of the truth tables of all bracketed formulae with $n$ distinct propositional variables
 connected by the binary connective of implication, and Kleene implication. This work is an extension of \cite{V1}, and \cite{V2}, 
 and also, a generalisation of a combinatorical structures that we have discussed in the former papers.
 For detailed enumerative and asymptotic results  please do refer back to \cite{V1} and \cite{V2}. Here we briefly discuss the basic definitions and results of our former results.
\\\\
We represent truth values of propositional variables and formulae
by $1$ (for ``true''),  $0$ (for ``false''), and $2$ (for ``unknown''). In this paper 
a propositional function of $n$ variables $p_1,\ldots,p_n$ is simply a function from $\{0,1,2\}^n$ to 
$\{0,1,2\}$. In two valued classical logic, it is well known that there are  $2^{2^n}$ propositional functions, 
and in three valued propositional logic, there are $3^{3^n}$ propositional functions, each of which can be represented by a
 formula involving the connectives $\neg$, $\vee$ and $\wedge$. A valuation is an assignment of values to the
variables $p_1,\ldots,p_n$, with consequent assignment of values to formulae.

The function represented by a formula is conveniently calculated using a truth
table. Each row of the truth table corresponds to a valuation.

We are interested in \emph{bracketed Kleene implications}, which are formulae
obtained from $p_1\Rightarrow p_2\Rightarrow\cdots\Rightarrow p_n$ by inserting brackets so that
the result is well-formed; the binary connective $\Rightarrow$ (``implies'') is
defined by the rule that, for any valuation $\nu$, and any propositional statement $\phi$ and $\psi$ 

\begin{displaymath}
\begin{array}{|c|c| c|c|}
\hline
\nu(\phi
\Rightarrow \psi) &1 & 0 & 2  \\ 
\hline 
 1 & 1 & 0 &2\\
\hline
0& 1 & 1&1\\
\hline
2 & 1 & 2 & 2\\
\hline

\end{array}
\end{displaymath}
\;\\
Here, for example, are the truth tables for the two bracketed implications
in $n=3$ variables.
\begin{figure}[h]

\centering

        \includegraphics[scale=0.7]{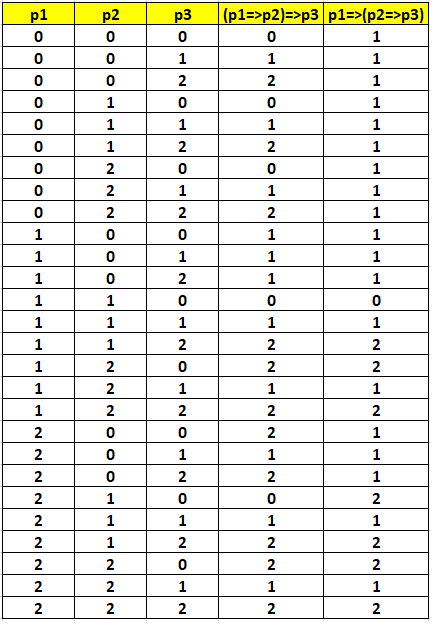}

 \caption{Kleeene implication, $n=3$}
\end{figure}

\begin{prop} \cite{V1} \;\;  Let $f_n,\; t_n$ and $\;u_n$ be the number of rows with the value “false” ,  “true”, and  “unknown”  in the Kleene truth tables of all
bracketed formulae with n distinct propositions $p_1, . . . , p_n$ connected by the binary connective
of implications. Then $u_n$, $f_n$, and $t_n$ have the following generating functions respectively:
\begin{equation} 
U(x)= \frac{1-\sqrt{1-12x}}{6}
\end{equation}

\begin{equation} 
F(x) = \frac{-2-\sqrt{1-12x}+\sqrt{5+24x+4\sqrt{1-12x}}}{6}
\end{equation}

\begin{equation}
T(x)=\frac{4-\sqrt{1-12x}-\sqrt{5+24x+4\sqrt{1-12x}}}{6}
\end{equation}

Where $G(x)$ is the generating function for the total number of entries  in all truth
tables for bracketed implications with $n$ distinct variables $p_1, \ldots ,p_n$, in Kleene's propositional calculus.
\begin{equation}G(x)= \frac{1-\sqrt{1-12x}}{2}
\end{equation}
\end{prop}

We can think of Kleene type truth tables connected by binary connective of implications  as a $C_n=\frac{1}{n}{{2n-2}\choose{n-1}}$ by $3^n$ rectangular array and count the number of truth, false and unknown entries. This gives us a sequence space which we can name it as $\mathcal{G}_3$; moreover $\mathcal{G}_2$ will be a subspace of $\mathcal{G}_3$ which is generated by $R(x)$ and $S(x)$ functions:
\begin{prop}\cite{V0}
Let $s_n$ and $\;r_n$ be the number of rows with the value “false” , and  “true” in the truth tables of all
bracketed formulae with $n$ distinct propositions $p_1, . . . , p_n$ connected by the binary connective
of implications. Then ${s}_n$,  and $r_n$ have the following generating functions respectively:
\begin{equation}
S(x)= \frac{-1-\sqrt{1-8x} + \sqrt{2+2\sqrt{1-8x}+8x}}{4}
\end{equation}
\begin{equation}
R(x)=\frac{3-\sqrt{1-8x}-\sqrt{2+2\sqrt{1-8x}+8x}}{4}
\end{equation}
\end{prop}
 
Beyond their analytic structures these sequence spaces also have an algebraic stucture. Namely they are commutative monoids. In this extension paper I will prove my claim. Of course, one needs to study immediately the Green Relations on such an algebraic structure and we will underline the GRs on these rectangular arrays. Here is an example of $C_4=5$ by $2^4$ a rectangular array, (from classical logic):\\\\
\begin{figure}[h]

\centering

        \includegraphics[scale=0.35]{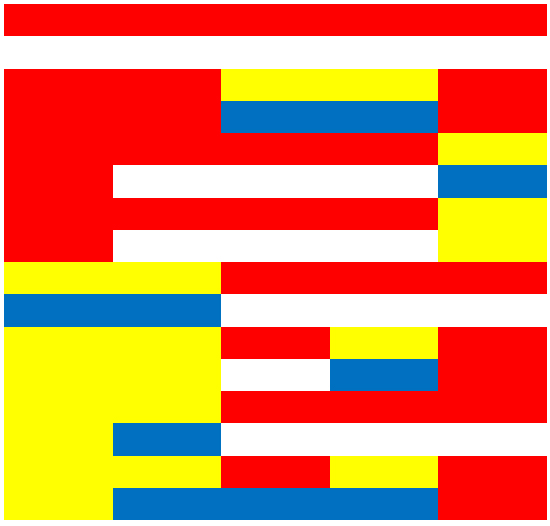}

 \caption{Binary connective of classical implication, $n=4$}

\end{figure}

Figure 2 is colour coded where the generating function is $R(x)R(x)$ counts the number of red entries; the generating function $R(x)S(x)$ counts the number of white entries; the generating function $S(x)R(x)$ counts of the number of yellow entries; and finally, the generating function $S(x)S(x)$ counts blue entries. This rectangular array is embedded in a bigger rectangular array of $\mathcal{G}_3$ at $n=4$, where the total number of entries equals to $3^4\times 5$. The factor 3 comes from additional the third truth value namely {\it{unknows}} .

\section{Algebraic structure}
Generating functions  $T(x)$, $F(x)$, and $U(x)$ form a commutative monoid with the multiplication. 
First we need to observe the following property from the following `partial' Cayley table
\begin{displaymath}
\begin{array}{|c||c|c| c|c|}
\hline
\times &I(x)&T(x) & F(x) & U(x)  \\
\hline 
\hline
I(x) &I(x) & T(x)& F(x)&U(x)\\
\hline
T(x)&T(x)& T(x)^2 & T(x)F(x)& T(x)U(x)\\
\hline
F(x)& F(x)& F(x)T(x) & F(x)^2 & F(x)U(x)\\
\hline
U(x)&U(x) & U(x)T(x) & U(x)F(x) & U(x)^2\\
\hline
\end{array}
\end{displaymath}
Where the identity generating function is defined as follows:
\[
I(x)=1 +0x+0x^2+0x^3+...
\]

\begin{deff}
\[
\bold{K}_3=\{ T(x), F(x), U(x), I(x) \}
\]
\end{deff}

\begin{thm}
Consider the following set 
\[
\mathcal{G}_3=\bigg\{ A(x): \;\;a_n< g_n,\;n> 1,\; \;A(x)\in\{B(x)^a\times C(x)^b: B(x),C(x)\in \bold{K}_3, \; a,b\geq 1\}\bigg\} 
\]
where $g_n=3^nC_n$, 
then $(\mathcal{G}_3, \times)$ is a commutative monoid.
\begin{enumerate}
\item \[ [x^n](A(x)B(x))=[x^n](B(x)A(x)),\;\;\; \forall A(x),B(x)\in \mathcal{G}_3\] 
\item \[ [x^n]A(x)^m < [x^n]G(x),\;\;\; \forall A(x)\in \mathcal{G}_3, \; \forall m\geq1 \] 
\item \[ [x^n](A(x)^kB(x)^l)< [x^n]G(x),\;\;\; \forall A(x),B(x)\in \mathcal{G}_3, \;\; \forall k,l\geq1 \] 
\item  \[ [x^n]((A(x)^kB(x)^l)C(x)^m)=[x^n](A(x)^k(B(x)^lC(x)^m) )< [x^n]G(x),\;\;\; \forall A(x),B(x),C(x)\in \mathcal{G}_3, \;\; \forall k,l,m\geq1 \] 
\end{enumerate}
\end{thm}
\begin{proof}Proving (1) and (4) are straightforward. Here we first prove (2), then (3).
The results follow by induction on $k>1$ and expressing the generating functions: 
\[
3U(x)^{k}= U(x)^{k-1}-xU(x)^{k-2} 
\]
\[
F(x)^{k}=2F(x)^{k-1}U(x)-F(x)^{k-1}+xF(x)^{k-2} 
\]
and
\[
T(x)^{k}=2/3 T(x)^{k-1}G(x)^2-2/3T(x)^{k-1} G(x)F(x)+T(x)^{k-1}F(x)^2 
\]
we have 
\[
[x^k]U(x)^n<[x^k]G(x),\;\;\;\; [x^k]F(x)^n<[x^k]G(x),\;\;\;\; [x^k]T(x)^n<G(x).
\]
Now we can prove (3). Consider $U(x)^mF(x)^k$, if $m>k$,  where $m,k\in\mathbb{N}$, then $U(x)^a(U(x)F(x))^k$ $\exists a \in\mathbb{N}$ and $a<k$.
Recall from our former paper, \cite{V2}, how we have defined $U(x)F(x)=U_2(x)$
\[
[x^n](U(x)^a(U(x)F(x))^k) = [x^n](U(x)^aU_2(x)^k) <[x^n]U(x)^k<[x^k]G(x).
\]
We can use the same trick for $U^m(x)T(x)^k$,
\[
[x^n](U(x)^mT(x)^m)=[x^n](U(x)^b(U(x)T(x))^k)=[x^n](U(x)^bT_5(x)^k)<[x^n]T(x)^k <[x^n]G(x)
\]
 and for $F^m(x)T^k(x)$:
\[
[x^n](F(x)^mT(x)^k) =[x^n](F^b(x)(F(x)T(x))^k) =[x^n](F(x)^bT_2(x)^k)<[x^n]T(x)^k<[x_n]G(x)
\]
\end{proof}
Therefore $T, F,$  and $\; U\;$  are well defined sets:
\begin{equation}\notag
\begin{split}
& T=\{ A(x): \; A(x)=T(x)^m, \; m\geq 1, \; a_n < g_n\}\\
& F=\{ B(x): \; B(x)=F(x)^m, \; m\geq 1, \; b_n < g_n\}\\
& U=\{ C(x): \; C(x)=U(x)^m, \; m\geq 1, \; c_n < g_n\}
\end{split}
\end{equation}
Then $\mathcal{K_T}=\bigg(T, I(x), \times \bigg)$, \; $\mathcal{K_F}=\bigg(F, I(x), \times \bigg)$,  $\;\mathcal{K_U}=\bigg(U, I(x), \times \bigg)$, and $\;\mathcal{K_I}=\bigg( I(x), \times \bigg)$
are commutative monoids.
Let $\mathcal{K}=\{T,F,U\}$
\[
\mathcal{G}_3=\prod_{\nu \in \mathcal{K}} \bold{K_\nu} 
\]
$\bigg(\mathcal{G}_3,I(x),\times\bigg)$is a commutative monoid.\\\\\\

Similarly in 2-valued classical logic,  please refer back to \cite{V2} and \cite{V0} for enumerative and asymptotic results,  $\mathcal{C_{R}}=\bigg(R, I(x), \times \bigg)$, \; $\mathcal{C_{S}}=\bigg(S, I(x), \times \bigg)$, and $\;\mathcal{C_{I}}=\bigg( I(x), \times \bigg)$
are commutative monoids in their own rights, and they are subset of $\mathcal{G}_3$ and hence submonoid of $\mathcal{G}_3$. 

\begin{equation}\notag
\begin{split}
& R=\{ A(x): \; A(x)=R(x)^m, \; m\geq 1, \; a_n < g_n\}\\
& S=\{ B(x): \; B(x)=S(x)^m, \; m\geq 1, \; b_n < g_n\}
\end{split}
\end{equation}

For classical logic we need to define $\mathcal{G}_2$: (here $g_n=2^nC_n$)
\[
\mathcal{G}_2=\bigg\{ A(x): \;\;a_n< g_n,\;n> 1,\; \;A(x)\in\{B(x)^a\times C(x)^b: B(x),C(x)\in \bold{K}_2, \; a,b\geq 1\}\bigg\} 
\]
Where
 \[
\bold{K}_2=\{ R(x), S(x), I(x) \}
\]
hence with a similar arguments $\bigg(\mathcal{G}_2,I(x),\times\bigg)$is a commutative monoid, and indeed a submonoid of $\mathcal{G}_3$.\\\\\\

Let $T$ be the underlying set of $\mathcal{K_T}$  then it can be shown with the above theorem, section (3), that $TA\subset \mathcal{K_T}$, where $A\in\mathcal{G}_3$. Thus $\mathcal{K_T}\mathcal{G}_3\subset \mathcal{K_T}$ and since multiplication is commutative we have a two sided ideal. Same is also true for other submonoids:
\[
\mathcal{K_T},\mathcal{K_F}, \mathcal{K_U}, \mathcal{K_I} \triangleleft \mathcal{G}_3.
\]
One can easily see that 
\[
\mathcal{G}_2\le \mathcal{G}_3\;\; \textit{ but }  \;\;\mathcal{G}_2\not\trianglelefteq \; \mathcal{G}_3
\]
and
\[
\mathcal{C_R}, \mathcal{C_S} \le \mathcal{G}_3, \;\; \textit{ but }  \;\; \mathcal{C_R},\mathcal{C_S} \not\trianglelefteq \; \mathcal{G}_3.
\]
However,
\[
\mathcal{C_R}, \mathcal{C_S} \le \mathcal{G}_2, \;\;\; \mathcal{C_R},\mathcal{C_S} \triangleleft \;\mathcal{G}_2.
\]\\
Consequently, we have the following result:\\\\\\
\begin{cor} 
Let $\mathcal{L, R}$ and $\mathcal{H}$ be the Green relations on $\mathcal{G}_3$ then $\mathcal{L=R=H}$.
\end{cor}
$\; $\\\\\\\\\\\\\\\\\\\\\\\\\\

I would like to thank to Prof F. Vivaldi for his valuable comments on the style and structure of the draft version of this paper.

\pagebreak

\;\\\\\\\\\\\\\\\\\\\\\\\\\\\\\\\\\\\\\\\\\\\\\\\\\\\\\\\\\\\\\\\\\\\\\\\\

\begin{flushright}
...Seni sevmek,\\
   Felsefedir kusursuz.\\
   İmandır, korkunç sabırlı.\\
   İp'in, kurşun'un rağmına,\\
   Yürür pervasız ve güzel.\\
   Sıradağları devirir,\\
   Akan suları çevirir,\\
   Alır yetimin hakkını,\\
   Buyurur, kitabınca...\\
A. Arif
\end{flushright}

\end{document}